\documentclass{article}
\usepackage{amsmath,amssymb,latexsym,amsthm,amscd}
\font\goth=eusm10

\newcommand\bc{\mathbf C}

\newcommand\E{\mathcal E}

\newcommand\bz{\mathbf{Z}}

\newcommand\bn{\mathbf{N}}

\newcommand\Oc{\hbox{\goth O}}

\newcommand\wH{\widetilde{H}}
\newcommand\wth{\widetilde{h}}

\newtheorem{theorem}{Theorem}

\newtheorem{corollary}{Corollary}

\newtheorem{lemma}{Lemma}

\numberwithin{proposition}{section}
\numberwithin{definition}{section}
\numberwithin{corollary}{section}
\numberwithin{remark}{section}
\numberwithin{lemma}{section}
\numberwithin{equation}{section}

\numberwithin{theorem}{section}

\numberwithin{question}{section}
\numberwithin{case}{section}
\numberwithin{example}{section}
\numberwithin{conjecture}{section}

\begin{document}
\title{On correspondences of a K3 surface\\ with itself. IV}
\author{C.G.Madonna \footnote{Supported by I3P contract} \ and
Viacheslav V.Nikulin \footnote{Supported by EPSRC grant
EP/D061997/1}}
\maketitle

\begin{abstract} Let $X$ be a K3 surface with a polarization $H$
of the degree $H^2=2rs$, $r,s\ge 1$, and the isotropic Mukai
vector $v=(r,H,s)$ is primitive. Moduli space of sheaves over $X$
with the isotropic Mukai vector $(r,H,s)$ is again a K3 surface,
$Y$. In \cite{Nik2} second author gave necessary and sufficient
conditions in terms of Picard lattice $N(X)$ of $X$ when $Y$ is
isomorphic to $X$ (important particular cases were considered in
\cite{Mad-Nik1}, \cite{Mad-Nik2} and \cite{Nik1}).

Here we show that these conditions imply existence of an
isomorphism between $Y$ and $X$ which is a composition of some
universal geometric isomorphisms between moduli of sheaves over
$X$, and geometric Tyurin's isomorphism between some moduli of
sheaves over $X$ and $X$ itself. It follows that for a general K3
surface $X$ with $\rho(X)=\text{rk\ }N(X)\le 2$ one has $Y\cong X$
if and only if there exists an isomorphism $Y\cong X$ which is a
composition of the universal and the Tyurin's isomorphisms.
\end{abstract}

\section{Introduction} \label{introduction}

Let $X$ be a K3 surface over $\bc$ with a polarization $H$ of the
degree $H^2=2rs$ where $r,s\in \bn$. Assume that the isotropic
Mukai vector $v=(r,H,s)$ is primitive. We denote by $Y=M_X(r,H,s)$
the moduli space of stable rank $r$ sheaves over $X$ with first
Chern class $H$ and Euler characteristic $r+s$ (i. e. with the
Mukai vector $v$). By results of Mukai \cite{Muk1}, \cite{Muk2},
$Y$ is again a K3 surface.

In (\cite{Nik2}, Theorem 4.4), second author gave some sufficient
and necessary conditions in terms of Picard lattice $N(X)$ of $X$
when $Y=M_X(r,H,s)$ is isomorphic to $X$. See Theorem
\ref{maintheorem} in Sect. \ref{section2} for the exact
formulation. Important particular cases of the theorem were also
obtained in \cite{Mad-Nik1}, \cite{Mad-Nik2} and \cite{Nik1}. The
sufficient part of these results used global Torelli Theorem for
K3 surfaces \cite{PShShaf} and was not effective.

Here we show that sufficient conditions of (\cite{Nik2}, Theorem
4.4), see Theorem \ref{maintheorem}, imply that there exists an
isomorphism $Y\cong X$ which is the composition of some universal
geometric isomorphisms between moduli of sheaves over $X$ (see
their exact definitions in Lemmas \ref{lemmadelta}, \ref{lemmaTD}
and Theorem \ref{maintheorem1}), and Tyurin's isomorphism between
some moduli of sheaves over $X$ and $X$ itself (see  Lemma
\ref{lemmaTyu}). The exact statement is given in Theorem
\ref{maintheorem2}.

Sufficient conditions of (\cite{Nik2}, Theorem 4.4), see Theorem
\ref{maintheorem}, are also necessary for a general K3 surface $X$
with Picard number $\rho(X)=\text{rk\ }N(X)\le 2$. Thus, the
preceding result implies that for a general K3 surface $X$ with
$\rho(X)\le 2$ one has that $Y\cong X$ if and only if there exists
an isomorphism $Y\cong X$ which is a composition of the universal
and Tyurin's isomorphisms. See Theorem \ref{maintheorem3} and
Corollary \ref{corollarymain} for the exact formulations.

These results show that the universal and Tyurin's isomorphisms
above are universally important in this problem (When is $Y\cong
X$?). They are sufficient to solve this problem for all general (for their
Picard lattice) K3 surfaces $X$ with $\rho(X)\le 2$.
Similar particular results were obtained in \cite{Mad-Nik3}
(geometric interpretation of \cite{Mad-Nik1}) and \cite{Mad-Nik4}
(geometric interpretation of \cite{Nik1}).


\section{Reminding of the main result of \cite{Nik2}}
\label{section2}

We denote by $X$ an algebraic K3 surface
over the field $\bc$ of complex numbers. I.e. $X$ is a non-singular
projective algebraic surface over $\bc$ with the trivial canonical class
$K_X = 0$ and the vanishing  irregularity $q(X)=0$.

We denote by $N(X)$
the Picard lattice (i.e. the lattice of 2-dimensional algebraic
cycles) of $X$. By $\rho (X)=\text{rk}\ N(X)$ we denote the Picard number of
$X$. By
\begin{equation}
T(X)=N(X)_{H^2(X,\bz)}^\perp
\end{equation}
we denote the
transcendental lattice of $X$.

\par\medskip

For a Mukai vector $v=(r,c_1,s)$ where $r\in \bn$, $s \in \bz$,
and $c_1\in N(X)$, we denote by $Y=M_X(r,c_1,s)$ the moduli space
of stable (with respect to some ample $H^\prime\in N(X)$) rank $r$
coherent sheaves on $X$ with first Chern class $c_1$, and Euler
characteristic $r+s$.

By results of Mukai \cite{Muk1}, \cite{Muk2}, under suitable
conditions on the Chern classes, the moduli space $Y$ is always
deformations equivalent to a Hilbert scheme of 0-dimensional
cycles on $X$ (of some dimension).
\par\medskip

The general common divisor of the Mukai vector $v$ is
$$
(r,c_1,s)=(r,d,s)
$$
if $c_1=dc_1^\prime$ where $d\in\bz$, $d\ge 0$ and $c_1^\prime$ is
primitive in the Picard lattice $N(X)$ which is a free
$\bz$-module. Here $c_1^\prime\in N(X)$ is {\it primitive} means
that $c_1^\prime/n\notin N(X)$ for any natural $n\ge 2$. A Mukai
vector $v$ is called {\it primitive} if its general common divisor
is one.

Let $X$ be a K3 surface with a polarization $H$ of the degree
$H^2=2rs$ where $r,s\in \bn$. Assume that the isotropic Mukai
vector $v=(r,H,s)$ is primitive. Then $Y=M_X(r,H,s)$ is another K3
surface. We are interested in the case when $Y\cong X$.

We denote $c=(r,s)$ and $a=r/c$, $b=s/c$. Then $(a,b)=1$. Let $H$
be divisible by $d\in \bn$ where $\wH=H/d$ is primitive in $N(X)$.
Primitivity of $v=(r,H,s)$ means that $(r,d,s)=(c,d)=1$. Since
$\wH^2=2abc^2/d^2$ is even, we have $d^2|abc^2$. Since
$(a,b)=(c,d)=1$, it follows that $d=d_ad_b$ where $d_a=(d,a)$,
$d_b=(d,b)$, and we can introduce integers
$$
a_1=\frac{a}{d_a^2},\ \ b_1=\frac{b}{d_b^2}\,.
$$
Then we obtain that $\wH^2=2a_1b_1c^2$.

Let $\gamma=\gamma (\wH)$ is defined by $\wH\cdot N(X)=\gamma
\bz$, i.e. $H\cdot N(X)=\gamma d\bz$. Clearly, $
\gamma|\wH^2=2a_1b_1c^2$.

We denote
\begin{equation}
\label{n(v)1} n(v)=(r,s,d\gamma)=(r,s,\gamma).
\end{equation}
By Mukai \cite{Muk1}---\cite{Muk3}, $T(X)\subset T(Y)$, and
$n(v)=[T(Y):T(X)]$ where $T(X)$ and $T(Y)$ are transcendental
lattices of $X$ and $Y$. Thus,
\begin{equation}
\label{n(v)2} Y\cong X\ \implies \
n(v)=(r,s,d\gamma)=(c,d\gamma)=(c,\gamma)=1.
\end{equation}

Assuming that $Y\cong X$ and then $n(v)=1$, we have
$\gamma|2a_1b_1$, and we can introduce
\begin{equation}
\gamma_a=(\gamma,a_1),\ \gamma_b=(\gamma,b_1),\
\gamma_2=\frac{\gamma}{\gamma_a\gamma_b}.
\end{equation}
Clearly, $\gamma_2|2$.

In (\cite{Nik2}, Theorem 4.4) the following general theorem had
been obtained (see its important particular cases in
\cite{Mad-Nik1}, \cite{Mad-Nik2} and \cite{Nik1}). In the theorem,
we use notations $c$, $a$, $b$, $d$, $d_a$, $d_b$, $a_1$, $b_1$
introduced above. The same notations $\gamma$, $\gamma_a$,
$\gamma_b$ and $\gamma_2$ as above are used if one replaces $N(X)$ by a
2-dimensional primitive sublattice
$N\subset N(X)$, e. g. $\wH\cdot N=\gamma \bz$,
$\gamma>0$. We denote $\det{N}=-\gamma\delta$ and $\bz f(\wH)$
denotes the orthogonal complement to $\wH$ in $N$.

\begin{theorem}
\label{maintheorem} Let $X$ be a K3 surface and $H$ a polarization
of $X$ such that $H^2=2rs$ where $r,s\in \bn$. Assume that the
Mukai vector $(r,H,s)$ is primitive. Let $Y=M_X(r,H,s)$ be the K3
surface which is the moduli of sheaves over $X$ with the isotropic
Mukai vector $v=(r,H,s)$. Let $\wH=H/d$, $d\in\bn$, be the
corresponding primitive polarization.

We have $Y\cong X$ if there exists $\wth_1 \in N(X)$ such that
$\wH$ and $\wth_1$ belong to a 2-dimensional primitive sublattice
$N \subset N(X)$ such that $\wH\cdot N=\gamma \bz$, $\gamma>0$,
$(c,d\gamma)=1$, and  the element $\wth_1$ belongs to the
$a$-series or to the $b$-series described below:

$\wth_1$ belongs to the $a$-series if
\begin{equation}
\wth_1^2=\pm 2b_1c, \ \ \wH\cdot \wth_1\equiv 0\mod
\gamma(b_1/\gamma_b)c,\ \ f(\wH)\cdot \wth_1\equiv 0\mod \delta
b_1c \label{aseries1}
\end{equation}
(where $\gamma_b=(\gamma,b_1)$);

$\wth_1$ belongs to the $b$-series if
\begin{equation}
\wth_1^2=\pm 2a_1c, \ \ \wH\cdot \wth_1\equiv 0\mod
\gamma(a_1/\gamma_a)c,\ \ f(\wH)\cdot \wth_1\equiv 0\mod\delta
a_1c \label{bseries1}
\end{equation}
(where $\gamma_a=(\gamma,a_1)$).

These conditions are necessary to have $Y\cong X$ if $\rho (X)\le
2$ and $X$ is a general K3 surface with its Picard lattice, i. e
the automorphism group of the transcendental periods
$(T(X),H^{2,0}(X))$ is $\pm 1$.
\end{theorem}

{\bf Remark 2.1.} In \cite{Nik2}, some additional primitivity
conditions on $\wth_1$ were required in \eqref{aseries1},
\eqref{bseries1}. From the proof of Theorem \ref{maintheorem2}
below, it will follow that they are unnecessary as sufficient
conditions. They are important as necessary conditions.

\medskip

The sufficient part of the proof of Theorem \ref{maintheorem} in
\cite{Nik2} (and its important particular cases in
\cite{Mad-Nik1}, \cite{Mad-Nik2} and \cite{Nik1}) used global
Torelli Theorem for K3 surfaces \cite{PShShaf}. I. e., under
conditions of Theorem \ref{maintheorem}, it was proved that the K3
surfaces $Y$ and $X$ have isomorphic periods. By global Torelli
Theorem for K3 surfaces \cite{PShShaf}, then $Y$ and $X$ are
isomorphic. In Sect. \ref{section3} below we will give a geometric
construction of the isomorphism between $Y$ and $X$. This will be
the main result of the paper.

In \cite{Mad-Nik3} and \cite{Mad-Nik4} similar geometric
construction was given in particular cases of Theorem
\ref{maintheorem} obtained in \cite{Mad-Nik1} and \cite{Nik1}
respectively.

\section{Geometric interpretation of Theorem \ref{maintheorem}}
\label{section3}

We use notations of the previous Section \ref{section2}.

We shall use the natural isomorphisms between moduli of sheaves
over a K3 surface $X$.

\begin{lemma}
\label{lemmadelta} Let $v=(r,H,s)$ be a primitive isotropic Mukai
vector for a K3 surface $X$, and $r,s\ge 1$. Then one has an
isomorphism, called {\bf reflection}, see \cite{Tyurin1},
\cite{Tyurin2},
$$
\delta:M_X(r,H,s)\cong M_X(s,H,r).
$$
 \end{lemma}

In (\cite{Tyurin1}, (4.11)) and (\cite{Tyurin2}, Lemma 3.4), the
geometric construction of the reflection $\delta$ is given under
the condition that $M_X(r,H,s)$ contains a regular bundle (it is
valid even if the Mukai vector $v=(r,H,s)$ is not isotropic). On
the other hand, by global Torelli Theorem for K3 surfaces
\cite{PShShaf}, existence of such isomorphism for isotropic
primitive $v$ is obvious. See similar proof of Theorem
\ref{maintheorem1} below.

\begin{lemma}
\label{lemmaTD}
Let $(r,H,s)$ be a Mukai vector for a K3 surface $X$ and $D\in N(X)$.
Then one has the natural isomorphism given by the tensor product
$$
T_D:M_X(r,H,s)\cong M_X(r,H+rD,s+r(D^2/2)+D\cdot H),\ \  \E\mapsto \E\otimes \Oc(D).
$$
Moreover, here Mukai vectors
$$
v=(r,H,s),\ \ v_1=(r,H+rD,s+r(D^2/2)+D\cdot H)
$$
have the same general common divisor and the same square under Mukai pairing.
In particular, they are primitive and isotropic simultaneously.
 \end{lemma}

\medskip

We also use the isomorphisms between moduli of sheaves over a K3
surface and the K3 surface itself found by A.N. Tyurin in
(\cite{Tyurin2}, Lemma 3.3).

\begin{lemma}
\label{lemmaTyu}
For a K3 surface $X$ and $h_1\in N(X)$ such that $\pm h_1^2>0$
and
\begin{equation}
h^0\Oc(h_1)=h^0\Oc(-h_1)=0\  if\  h_1^2<0,
\label{Tyurin1}
\end{equation}
there is a geometric Tyurin's isomorphism
$$
Tyu(\pm h_1):M_X(\pm h_1^2/2,h_1, \pm 1)\cong X .
$$
\end{lemma}

Like in Theorem \ref{maintheorem1} below, using global Torelli
Theorem for K3 surfaces \cite{PShShaf}, one can show that even
without Tyurin's condition \eqref{Tyurin1} always there exists
some isomorphism
$$
M_X(\pm h_1^2/2,h_1, \pm 1)\cong X.
$$
We call such isomorphisms as {\it Tyurin's isomorphisms} either.
When $h_1$ also satisfies \eqref{Tyurin1}, there exists a direct
Tyurin's geometric construction of some isomorphism $M_X(\pm
h_1^2/2,h_1,\pm 1)\cong X$.

\medskip

We shall use the following result which had been proved implicitly
in \cite{Nik2}. See Sect. 2.3 and the proof of Theorem 2.3.3 in
\cite{Nik2}. In \cite{Nik2} much more difficult results related to
arbitrary Picard lattice had been considered, and respectively the
proofs were long and difficult. Therefore, below we also give a
much simpler proof of this result.

\begin{theorem}
Let $v=(r,H,s)$  be an isotropic Mukai vector on a K3 surface $X$
where $r,\ s \in \bn$, $H\in N(X)$, $H^2=2rs$,  and $H$ is
primitive (then $v$ is also primitive).

Let $d_1,d_2\in \bn$ and
$(d_1,s)=(d_2,r)=(d_1,d_2)=1$.

Then the Mukai vector
$v_1=(d_1^2r, d_1d_2H,d_2^2s)$
is also primitive and isotropic,
and  there exists a natural isomorphism of moduli of sheaves
$$
\nu(d_1,d_2): M_X(r,H,s)\cong M_X(d_1^2r,d_1d_2H,d_2^2s).
$$
\label{maintheorem1}
\end{theorem}

\begin{proof} Let us consider the case $(d_1,d_2)=(d,1)$ where $(d,s)=1$ (general case
is similar).

Let $c=(r,s)$ and $a=r/c$, $b=s/c$. Then $v=(ac,H,bc)$ where
$H^2=2abc^2$ and $H$ is primitive. By global Torelli Theorem
\cite{PShShaf}, it is enough to show that periods of
$Y=M_X(r,H,s)$ and $Y_1=M_X(d^2r,dH,s)$ are isomorphic.

By results of Mukai \cite{Muk1}, \cite{Muk2}, cohomology of $Y$
(and similarly of $Y_1$) are equal to
$$
H^2(Y,\bz)=v^\perp/\bz v
$$
where we consider $v$ as the element of Mukai lattice
$H^\ast(X,\bz)$, and $H^{2,0}(Y)$ is the projection of
$H^{2,0}(X)$.

Using the variant of Witt's Theorem from \cite{PShShaf}, we can model the necessary
calculations as follows.

Let $U^{(1)}$ be an even unimodular hyperbolic plane  with the
basis $e_1,\ e_2$ where $e_1^2=e_2^2=0$ and $e_1\cdot e_2=-1$. Let
$U^{(2)}$ be another even unimodular hyperbolic plane with the
bases $f_1,\ f_2$ where $f_1^2=f_2^2=0$ and $f_1\cdot f_2=1$. We
consider the orthogonal sum $U^{(1)}\oplus U^{(2)}$ (the model of
$H^\ast (X,\bz)$; to get $H^\ast(X,\bz)$, one should add to
$U^{(1)}\oplus U^{(2)}$ an unimodular even lattice of signature
$(2, 18)$ which is the same for all three $X$, $Y$ and $Y_1$).
Then
$$
v=ace_1+bce_2+H,\ \ H=abc^2f_1+f_2,
$$
thus
$$
v=ace_1+bce_2+abc^2f_1+f_2.
$$
Then $N(X)=\bz H$ models the Picard lattice of $X$ and $T(X)=\bz t$, $t=-abc^2f_1+f_2$
models the transcendental lattice of $X$.

We have $\xi=x\,e_1+y\,e_2+z\,f_1+w\,f_2\in v^\perp$ if and only
if $-bcx-acy+z+abc^2w=0$, equivalently $z=bcx+acy-abc^2w$ where
$x,y,w,z\in \bz$. Thus,
$$
\xi=x(e_1+bc\,f_1)+y(e_2+ac\,f_1)+w(-abc^2\,f_1+f_2)\,
$$
and $\alpha=e_1+bc\,f_1$, $\beta=e_2+ac\,f_1$, $t=-abc^2\,f_1+f_2$
give the basis of $v^\perp$. We have $v=ac\,\alpha+bc\,\beta+t$
and
$$
ac\,\alpha\hskip-7pt\mod {\bz v}+bc\,\beta\hskip-7pt\mod{\bz v}+
t\hskip-7pt\mod{\bz v}=0.
$$
It follows that $\overline{\alpha}=\alpha\hskip-5pt\mod {\bz v}$,
$\overline{\beta}= \beta\hskip-5pt\mod{\bz v}$ give a basis of
$v^\perp/\bz v$ which models $H^2(Y,\bz)$. We have
$\overline{\alpha}^2=\overline{\beta}^2=0$ and
$\overline{\alpha}\cdot \overline{\beta}=-1$. Thus $v^\perp/\bz
v\cong U$. We see that $\overline{t}=t\hskip-5pt\mod{\bz
v}=-ac\overline{\alpha}-bc\overline{\beta}$ and then
$\widetilde{t}=\overline{t}/c=-a\overline{\alpha}-b\overline{\beta}\in
(v^\perp/\bz v)$, and $\bz \widetilde{t}$ models the
transcendental lattice $T(Y)$ of $Y$. Its orthogonal complement
$\bz h$ where $h=a\overline{\alpha}-b\overline{\beta}$ models the
Picard lattice $N(Y)$ of $Y$. We have $h^2=2ab >0$.

Let us make similar calculations for $Y_1$. We have
$$
v_1=d^2ac\,e_1+bc\,e_2+dH,\ \ H=abc^2\,f_1+f_2,
$$
and
$$
v_1=d^2ac\,e_1+bc\,e_2+dabc^2\,f_1+d\,f_2.
$$
We have $\xi_1=x_1e_1+y_1e_2+z_1f_1+w_1f_2\in v_1^\perp$ if and
only if $-bcx_1-d^2acy_1+dz_1+dabc^2w_1=0$. Since $(d,bc)=1$, we
obtain $x_1=d\widetilde{x}_1$, and
$z_1=bc\widetilde{x}_1+dacy_1-abc^2w_1$ where
$\widetilde{x}_1,y_1,w_1,z_1\in \bz$. Then
$$
\xi_1=\widetilde{x}_1(de_1+bcf_1)+y_1(e_2+dacf_1)+w_1(-abc^2f_1+f_2)
$$
and $\alpha_1=d\,e_1+bc\,f_1$, $\beta_1=e_2+dac\,f_1$,
$t=-abc^2\,f_1+f_2$ give the basis of $v_1^\perp$. We have
$v_1=dac\,\alpha_1+bc\,\beta_1+d\,t$, and
$$
dac\,\alpha_1\hskip-7pt\mod {\bz v_1}+bc\,\beta_1\hskip-7pt\mod{\bz v_1}+
d\,t\hskip-7pt\mod{\bz v_1}=0.
$$
Since $(d,bc)=1$, we see that $t\hskip-5pt\mod{\bz v_1}=c\tilde{t}_1$,
$\beta_1\hskip-5pt\mod{\bz v_1}=d\tilde{\beta}_1$ where
$\tilde{t}_1, \tilde{\beta}_1\in v_1^\perp/\bz v_1$, and
$$
a\,\alpha_1\hskip-7pt\mod {\bz v_1}+b\tilde{\beta}_1+\tilde{t}_1=0.
$$
We have  $(\alpha_1\hskip-5pt\mod {\bz v_1})^2=\tilde{\beta}_1^2=0$ and
$(\alpha_1\hskip-5pt\mod {\bz v_1})\cdot \tilde{\beta}_1=-1$.
Thus, again $\alpha_1\hskip-5pt\mod {\bz v_1}$, $\tilde{\beta}_1$ give canonical
generators of the unimodular lattice $U$.
Then they give a basis of $v_1^\perp/\bz v_1$ which models $H^2(Y_1)$. Moreover
$\bz \tilde{t}_1$ where
$\tilde{t}_1= -a\,\alpha_1\hskip-7pt\mod {\bz v}-b\,\tilde{\beta}_1$ (from above),
 models the transcendental lattice of $Y_1$.
Its orthogonal complement $\bz h_1$ where
$h_1=a\,\alpha_1\hskip-5pt\mod {\bz v}-b\,\tilde{\beta}_1$ models the
Picard lattice of $Y_1$.

We see that our descriptions for $Y$ and $Y_1$ above are evidently
isomorphic if we identify $\overline{\alpha}$, $\overline{\beta}$
with $\alpha_1\hskip-5pt\mod {\bz v}$ and $\tilde{\beta}_1$
respectively. This shows that $Y$ and $Y_1$ have isomorphic
periods and are isomorphic by global Torelli Theorem for K3
surfaces \cite{PShShaf}.

This finishes the proof.
\end{proof}

\medskip

{\bf Remark 3.1.} It would be very interesting to find a direct
geometric proof of Theorem \ref{maintheorem1} which does not use
global Torelli Theorem for K3 surfaces. It seems, considerations
by Mukai in  \cite{Muk3} are related to this problem (especially
see Theorem 1.2  in  \cite{Muk3}).

On the other hand, the isomorphism $\nu(d_1,d_2)$ is very universal, and
it exists even for general (with Picard number one) K3 surfaces.
By \cite{Nik0}, there exists only one isomorphism (or two isomorphisms
for the degree two) between algebraic K3 surfaces with Picard number one.
Thus, we can consider this isomorphism as geometric by definition.

\medskip

We will show below that under the conditions of Theorem
\ref{maintheorem}, there exists an isomorphism between
$Y=M_X(r,H,s)$ and $X$ which is the composition of the universal
geometric isomorphisms above. We use notations of Section
\ref{section2} and Theorem \ref{maintheorem}.

\begin{theorem}
\label{maintheorem2} Let $X$ be a K3 surface and $H$ a
polarization of $X$ such that $H^2=2rs$ where $r,s\in \bn$. Assume
that the Mukai vector $(r,H,s)$ is primitive. Let $Y=M_X(r,H,s)$
be the K3 surface which is the moduli of sheaves over $X$ with the
isotropic Mukai vector $v=(r,H,s)$. Let $\wH=H/d$, $d\in \bn$, be
the corresponding primitive polarization.

Assume that there exists $\wth_1 \in N(X)$ such that $\wH$ and
$\wth_1$ belong to a 2-dimensional primitive sublattice $N \subset
N(X)$ such that $\wH\cdot N=\gamma \bz$, $\gamma>0$,
$(c,d\gamma)=1$, and the element $\wth_1$ belongs to the
$a$-series or to the $b$-series described below:

$\wth_1$ belongs to the $a$-series if
\begin{equation}
\wth_1^2=\pm 2b_1c, \ \ \wH\cdot \wth_1\equiv 0\mod
\gamma(b_1/\gamma_b)c,\ \ f(\wH)\cdot \wth_1\equiv 0\mod \delta
b_1c \label{aseries2}
\end{equation}
(where $\gamma_b=(\gamma,b_1)$);

$\wth_1$ belongs to the $b$-series if
\begin{equation}
\wth_1^2=\pm 2a_1c, \ \ \wH\cdot \wth_1\equiv 0\mod
\gamma(a_1/\gamma_a)c,\ \ f(\wH)\cdot \wth_1\equiv 0\mod\delta
a_1c \label{bseries2}
\end{equation}
(where $\gamma_a=(\gamma,a_1)$).

Then we have:

If $\wth_1$ belongs to the $a$-series, then
\begin{equation}
\wth_1=d_2\wH+b_1c D\ \text{for some\ } d_2\in \bn,\  D\in N,
\label{aseriesnew}
\end{equation}
which defines the isomorphism
\begin{equation}
Tyu(\pm \wth_1)\cdot T_D\cdot \nu(1,d_2)\cdot \delta\cdot
\nu(d_a,d_b)^{-1}: Y=M_X(r,H,s)\cong X . \label{aseriesisomnew}
\end{equation}

If $\wth_1$ belongs to the $b$-series, then
\begin{equation}
\wth_1=d_2\wH+a_1cD\ \text{for some\ } d_2\in \bn,\ D\in N,
\label{bseriesnew}
\end{equation}
which defines the isomorphism
\begin{equation}
Tyu(\pm \wth_1)\cdot T_D\cdot \nu(1,d_2)\cdot \nu(d_a,d_b)^{-1}:
Y=M_X(r,H,s)\cong X . \label{bseriesisomnew}
\end{equation}
\end{theorem}

\begin{proof} We use the description of the pair $\wH\in N$ given in Proposition
3.1.1 in \cite{Nik2} (where one should replace $N(X)$ by $N$)
which directly follows from $\text{rk\ } N=2$, $\wH$ is primitive
in $N$,  and $\wH\cdot N=\gamma \bz$. We denote $\det
N=-\delta\gamma$, $\delta\in \bn$, and $\bz f(\wH)$ the orthogonal
complement to $\wH$ in $N$. Then
$f(\wH)^2=-2a_1b_1c^2\delta/\gamma$.

For some $\mu\in (\bz/(2a_1b_1c^2/\gamma))^\ast$ ($\pm \mu$ is the
invariant of the pair $\wH\in N$) where $\delta\equiv
\mu^2\gamma\mod {4a_1b_1c^2/\gamma}$, one has
\begin{equation}
N=[\wH,f(\wH), w=\frac{\mu \wH+f(\wH )}{2a_1b_1c^2/\gamma}],\
\label{latticeN1}
\end{equation}
where $[\ \ \cdot\  \ ]$ means `generated by $\cdot$ ".

It follows that
\begin{equation}
N=\{z=\frac{x\wH+yf(\wH)}{2a_1b_1c^2/\gamma }\ |\ x,y\in \bz\
\text{and}\ x\equiv \mu y\mod \frac{2a_1b_1c^2}{\gamma}\}
\label{latticeN2}
\end{equation}
where
\begin{equation}
z^2=\frac{\gamma x^2-\delta y^2}{2a_1b_1c^2/\gamma }\ .
\label{latticeN3}
\end{equation}

In the conditions of Theorem \ref{maintheorem1}, let us assume
that $\wth_1$ belongs to the $a$-series. Then
$\wth_1=(p\wH+qf(\wH))/(2a_1b_1c^2/\gamma )$ where $p,q\in \bz$
and
\begin{equation}
p\equiv \mu q\mod{\frac{2a_1b_1c^2}{\gamma}}\,. \label{latticeN4}
\end{equation}
By the condition of $a$-series, $\wH\cdot \wth_1=\gamma p\equiv
0\mod{\gamma (b_1/\gamma_b)c}$. Equivalently, $p\equiv 0\mod
{(b_1/\gamma_b)c}$. By \eqref{latticeN4}, then $q\equiv 0\mod
{(b_1/\gamma_b)c}$. Denoting $p=p_1(b_1/\gamma_b)c$ and
$q=q_1(b_1/\gamma_b)c$ where $p_1,\ q_1\in \bz$, we get that
$$
\wth_1=\frac{p_1\wH+q_1f(\wH)}{(2/\gamma_2)(a_1/\gamma_a)c}
$$
where $p_1\equiv \mu q_1\mod{(2/\gamma_2)(a_1/\gamma_a)c}$.

We have $f(\wH)=(2a_1b_1c^2/\gamma)\,w-\mu \wH$ where $w\in N$.
Then
\begin{equation}
\wth_1=\frac{p_1\wH+q_1((2a_1b_1c^2/\gamma) w-\mu
\wH)}{(2/\gamma_2)(a_1/\gamma_a)c}=\frac{p_1-\mu
q_1}{(2/\gamma_2)(a_1/\gamma_a)c} \wH+q_1(b_1/\gamma_b)c\, w\,.
\label{formula0}
\end{equation}
By the condition of $a$-series,
$$
f(\wH)\cdot \wth_1=-q_1(b_1/\gamma_b)\delta\,c \equiv 0\mod
{\delta b_1c},
$$
and we then get
$$
q_1\equiv 0\mod \gamma_b.
$$

Since $\mu$ is defined $\mod
(2/\gamma_2)(a_1/\gamma_a)(b_1/\gamma_b)c^2$ and $q_1\equiv 0\mod
\gamma_b$, then $p_1-\mu q_1\mod (2/\gamma_2)(a_1/\gamma_a)b_1c^2$
is correctly defined. Since $p_1-\mu q_1\equiv 0\mod
(2/\gamma_2)(a_1/\gamma_a)c$, then
$$
\frac{p_1-\mu q_1}{(2/\gamma_2)(a_1/\gamma_a)c}\mod b_1c
$$
is correctly defined. Since $q_1\equiv 0\mod \gamma_b$, we obtain
from \eqref{formula0} a correct congruence
\begin{equation}
\wth_1\equiv \frac{p_1-\mu
q_1}{(2/\gamma_2)(a_1/\gamma_a)c}\wH\mod {b_1c\,N}\, .
\label{formula1}
\end{equation}
There exists $d_2\in \bn$ such that
\begin{equation}
d_2\equiv \frac{p_1-\mu q_1}{(2/\gamma_2)(a_1/\gamma_a)c}\mod
b_1c\,. \label{formula2}
\end{equation}
Then
\begin{equation}
\wth_1=d_2\wH+b_1c D, \ d_2\in \bn,\  D\in N.
\label{formula3}
\end{equation}
By Theorem \ref{maintheorem1}, we have an isomorphism
$$
\nu(d_a,d_b)^{-1}:Y=M_X(ac,H,bc)\cong M_X(a_1c,\wH,b_1c)
$$
(this is one of the first steps of the proof of Theorem
\ref{maintheorem} in \cite{Nik2}).   By Lemma \ref{lemmadelta}, we
have an isomorphism
$$
\delta:M_X(a_1c,\wH,b_1c)\to M_X(b_1c,\wH,a_1c).
$$
By Lemma \ref{lemmaTD} and \eqref{formula3}, we obtain
$$
T_D:(b_1c,d_2\wH,d_2^2a_1c )\to (b_1c, \wth_1,\pm 1)
$$
since $(b_1c,d_2\wH,d_2^2a_1c)$ is isotropic Mukai vector and
$\wth_1^2=\pm 2b_1c$. Since $(b_1c,\wth_1,\pm 1)$ is evidently
primitive, $(b_1c,d_2\wH,d_2^2a_1c)$ is also primitive, and then
$(d_2,b_1c)=1$. By Theorem \ref{maintheorem1}, then we have an
isomorphism
$$
\nu (1,d_2):M_X(b_1c,\wH,a_1c)\cong M_X(b_1c,d_2\wH,d_2^2a_1c),
$$
and by Lemma \ref{lemmaTD}, an isomorphism
$$
T_D:M_X(b_1c,d_2\wH,d_2^2a_1c)\cong M_X(b_1c,\wth_1,\pm 1).
$$
At last, by Lemma \ref{lemmaTyu}, we have an isomorphism
$$
Tyu(\pm \wth_1):M_X(b_1c,\wth_1,\pm 1)\cong X.
$$
Thus, we obtain the isomorphism
\begin{equation}
Tyu(\pm \wth_1)\cdot T_D\cdot \nu(1,d_2)\cdot \delta\cdot
\nu(d_a,d_b)^{-1}: Y=M_X(r,H,s)\cong X . \label{isoma}
\end{equation}

If $\wth_1$ belongs to the $b$-series, similarly we obtain
\begin{equation}
\wth_1=d_2\wH+a_1cD,\ d_2\in \bn,\ D\in N,
\label{formula4}
\end{equation}
and we get the isomorphism
\begin{equation}
Tyu(\pm \wth_1)\cdot T_D\cdot \nu(1,d_2)\cdot \nu(d_a,d_b)^{-1}:
Y=M_X(r,H,s)\cong X . \label{isomb}
\end{equation}
We don't need the reflection $\delta$ in this case.

This finishes the proof of Theorem \ref{maintheorem2}.
\end{proof}

{\bf Remark 3.2.} Since conditions of Theorem \ref{maintheorem2}
have given an isomorphism $Y\cong X$, all other conditions of
Theorem 4.4 in \cite{Nik2} are unnecessary as sufficient
conditions. They are important as necessary conditions.

\medskip

 Since conditions of Theorem \ref{maintheorem2}
are also necessary for a general K3 surface $X$ with $\rho(X)\le
2$ to have $Y\cong X$ (see Theorem \ref{maintheorem}), we also
obtain a very interesting

\begin{theorem}
\label{maintheorem3}
Let $X$ be a K3 surface with a polarization
$H$ such that $H^2=2rs$, $r,s\ge 1$, the Mukai vector $(r,H,s)$ be
primitive, and $Y=M_X(r,H,s)$ be the moduli of sheaves over $X$
with the isotropic  Mukai vector $(r,H,s)$. Assume that
$\rho(X)\le 2$ and $X$ is general with its Picard lattice (i. e.
the automorphism group of the transcendental periods $Aut(T(X),
H^{2,0}(X))=\pm 1$). Let $\wH=H/d$, $d\in \bn$, be the corresponding
primitive polarization.

Then $Y=M_X(r,H,s)$ is isomorphic to $X$ if and only if there
exists $d_2\in \bn$ and $D\in N(X)$ such that

either
\begin{equation}
\wth_1=d_2\wH+b_1c D\ \text{has}\ \wth_1^2=\pm 2b_1c,\
\label{aseriesnew2}
\end{equation}
which defines the isomorphism
\begin{equation}
Tyu(\pm \wth_1)\cdot T_D\cdot \nu(1,d_2)\cdot \delta\cdot
\nu(d_a,d_b)^{-1}: Y=M_X(r,H,s)\cong X,\label{aseriesisomnew2}
\end{equation}

or
\begin{equation}
\wth_1=d_2\wH+a_1cD\ \text{has}\ \wth_1^2=\pm 2a_1c,
\label{bseriesnew2}
\end{equation}
which defines the isomorphism
\begin{equation}
Tyu(\pm \wth_1)\cdot T_D\cdot \nu(1,d_2)\cdot \nu(d_a,d_b)^{-1}:
Y=M_X(r,H,s)\cong X . \label{bseriesisomnew2}
\end{equation}
\end{theorem}

Theorem \ref{maintheorem3} gives an isomorphism between
$Y$ and $X$ as an exact composition of
the universal isomorphisms $\delta$, $T_D$, $\nu(d_1,d_2)$ and
$Tyu$.
Forgetting about the exact form of the composition, evidently we also obtain the following corollary.

\begin{corollary}
Let $X$ be a K3 surface with a polarization $H$ such that
$H^2=2rs$, $r,s\ge 1$, the Mukai vector $(r,H,s)$ be primitive,
and $Y=M_X(r,H,s)$ be the moduli of sheaves over $X$ with the
isotropic  Mukai vector $(r,H,s)$. Assume that $\rho(X)\le 2$ and
$X$ is general with its Picard lattice (i. e. the automorphism
group of the transcendental periods $Aut(T(X), H^{2,0}(X))=\pm
1$).

Then $Y=M_X(r,H,s)$ is isomorphic to $X$ if and only if there
exists an isomorphism between $Y=M_X(r,H,s)$ and $X$ which is a
composition of the universal isomorphisms $\delta$, $\nu(d_1,d_2)$
and $T_D$ between moduli of sheaves over $X$, and the universal
Tyurin's isomorphism $Tyu$ between a moduli of sheaves over $X$
and $X$ itself. \label{corollarymain}
\end{corollary}

Here it is important that the isomorphisms $\delta$, $T_D$ and
$Tyu$ have a geometric description which does not use Global
Torelli Theorem for K3 surfaces. Only for the isomorphism
$\nu(d_1,d_2)$ we don't know a geometric construction. On the
other hand, the isomorphism $\nu(d_1,d_2)$ is very universal and
is geometric by definition. See our considerations at the
beginning of this Section.

These results also show that Tyurin's isomorphisms $Tyu$ and
$\delta$, $T_D$, $\nu(d_1,d_2)$ as well, are universally important
in the problem (When is $Y=M_X(r,H,s)\cong X$?) which we had
considered. They are sufficient to solve this problem for all general (for their
Picard lattice) K3 surfaces $X$ with $\rho(X)\le 2$.

\

\

C.G.Madonna \par
Math. Dept.,
CSIC, C/ Serrano 121,
28006 Madrid,
SPAIN

carlo@madonna.rm.it \ \
cgm@imaff.cfmac.csic.es

\

\

V.V. Nikulin \par
Deptm. of Pure Mathem. The University of Liverpool, Liverpool\par
L69 3BX, UK;
\vskip1pt
Steklov Mathematical Institute,\par
ul. Gubkina 8, Moscow 117966, GSP-1, Russia

vnikulin@liv.ac.uk \ \
vvnikulin@list.ru

\end{document}